\documentclass{amsart}
\usepackage{graphicx}
\usepackage{amssymb}
\usepackage{epstopdf}
\usepackage{nicefrac}
\usepackage{enumerate}
\usepackage{hyperref}
\usepackage{float}
\usepackage{tabularx}
\usepackage[usenames,dvipsnames,svgnames,table]{xcolor}
 
\input xy
\xyoption {all}

\newtheorem{thm}{THEOREM}[section]

\newtheorem{cor}[thm]{COROLLARY}
\newtheorem{defn}[thm]{DEFINITION}

\newtheorem{lemma}[thm]{LEMMA}

\newtheorem{prop}[thm]{PROPOSITION}

\newcommand{\ds}{\displaystyle}

\newcommand{\F}{{\mathcal F}} 

\newcommand{\cG}{{\mathcal G}}
\newcommand{\cGF}{\cG_{\F}} 

\newcommand{\cI}{{\mathcal I}}
\newcommand{\cJ}{{\mathcal J}}

\newcommand{\cL}{{\mathcal L}}
\newcommand{\cM}{{\mathcal M}}

\newcommand{\cP}{{\mathcal P}}
\newcommand{\cR}{{\mathcal R}}
\newcommand{\cS}{{\mathcal S}}
\newcommand{\cT}{{\mathcal T}}
\newcommand{\cU}{{\mathcal U}}
\newcommand{\cV}{{\mathcal V}}

\newcommand{\dF}{d_{\F}} 

\newcommand{\dM}{d_{\fM}} %
\newcommand{\dTi}{d_{\cT_i}} %
\newcommand{\dX}{d_{\fX}} %
\newcommand{\e}{{\epsilon}} 
\newcommand{\eM}{{\epsilon_{\fM}}}

\newcommand{\fD}{{\mathfrak{D}}}

\newcommand{\fM}{{\mathfrak{M}}}
\newcommand{\fR}{{\mathfrak{R}}}
\newcommand{\fT}{{\mathfrak{T}}}

\newcommand{\fX}{{\mathfrak{X}}}
  
\newcommand{\lF}{{\lambda_{\mathcal F}}}

\newcommand{\mR}{{\mathbb R}}
\newcommand{\mS}{{\mathbb S}}

\newcommand{\mZ}{{\mathbb Z}}
\newcommand{\oU}{{\overline{U}}}
\newcommand{\vp}{{\varphi}}

\newcommand{\whU}{{\widehat U}}
\newcommand{\whvarp}{{\widehat \varphi}}

\newcommand{\psg}{{\rm pseudo}{\star}{\rm group}}

\newcommand{\eF}{{\epsilon_{\F}}}    
\newcommand{\eU}{{\epsilon_{\cU}}} 

\begin{document}


\title{Manifold--like matchbox manifolds}


\author{Alex Clark}
\address{Alex Clark, Department of Mathematics, University of Leicester, University Road, Leicester LE1 7RH, United Kingdom}
\email{Alex.Clark@leicester.ac.uk}

\author{Steven Hurder}
\address{Steven Hurder, Department of Mathematics, University of Illinois at Chicago, 322 SEO (m/c 249), 851 S. Morgan Street, Chicago, IL 60607-7045}
\email{hurder@uic.edu}

\author{Olga Lukina}
\address{Olga Lukina, Department of Mathematics, University of Illinois at Chicago, 322 SEO (m/c 249), 851 S. Morgan Street, Chicago, IL 60607-7045}
\email{lukina@uic.edu}

 \thanks{Version date: April 22, 2017}

\date{}


  \begin{abstract}
  A matchbox manifold is a generalized lamination, which is a continuum whose arc-components define the leaves of a foliation of the space.
The main result of this paper    implies that a matchbox manifold which is manifold-like must be homeomorphic to a weak solenoid. 
  \end{abstract}

\maketitle


\section{Introduction} \label{sec-intro}

A \emph{continuum} is a compact, connected, and non-empty
metrizable space. The notion of a manifold-like continuum is derived from the notion of an $\e$-map, which was introduced by 
   Alexandroff     \cite{Alexandroff1928}: 
\begin{defn}    \label{def-Ylike}
 Let  $X$ be a metric space, $Y$ a topological space and $\e > 0$ a constant. Then a   map $f \colon X \to Y$ is  said to be an $\e$-map if $f$ is a continuous surjection and for each point $y\in Y$,   the inverse image $f^{-1}(y)$   has diameter less than $\e$.   A metric space $X$ is said to be \emph{$Y$--like}, for some topological space $Y$, if for every   $\e>0$, there is an $\e$-map $f_\e \colon X \to Y$.
 \end{defn}

For example, a space $X$ is \emph{circle-like} if it is $Y$-like, where $Y = \mS^1$ is the circle. 
More generally, let 
\begin{equation}
\cM(n)=\{M\,|\, M \text{ is a closed connected manifold of dimension }n\,\}.
\end{equation}
\begin{defn}  \label{def-Mlike}
 A continuum  $X$ is said to be   \emph{manifold-like}, if there exists $n \geq 1$ such that for every   $\e>0$, there exists $M_{\e} \in \cM(n)$ and an $\e$-map   $f_\e \colon X \to M_{\e}$.
 \end{defn}

The study of the properties of $\e$-maps and   manifold-like continua has a long history in the study  of the topology of spaces.   Eilenberg showed in \cite{Eilenberg1938} that an $\e$-map, for $\e > 0$ sufficiently small,   admits a left approximate inverse.
 Ganea  studied the properties of  compact, locally connected  manifold-like ANR's of dimension $n$ in \cite{Ganea1959}, and showed that such a space  has the homotopy type of a closed $n$-manifold. 
Deleanu \cite{Deleanu1962,Deleanu1963} showed that an $n$-dimensional connected polyhedron which is manifold-like is a closed pseudo-manifold.
Bob Edwards gave in his 1978  ICM address \cite{Edwards1980}  an overview of the further applications of $\e$-approximations and homeomorphisms.

Marde\v{s}i\'{c} and Segal    \cite{MardesicSegal1963,MardesicSegal1967} studied the properties of manifold-like connected polyhedron, and gave conditions under which such spaces must be a topological manifold. These   authors used a technique of approximation of the given continuum by an inverse limit of spaces, and noted that their results do not apply to a continuum which is not locally connected, such as the   dyadic solenoid. 

 The goal of this work is to characterize a class of manifold-like continua for which the Marde\v{s}i\'{c} and Segal results do not apply. These are the matchbox manifolds as studied by the authors in 
 \cite{ClarkHurder2013,CHL2014,CHL2013c}, and discussed below. Our study of matchbox manifolds in these works was inspired by a result of Bing in \cite{Bing1960}.
 Recall that  a topological  space $X$ is \emph{homogeneous} if for
every $x, y \in X$, there exists a homeomorphism $h \colon X \to X$ such that $h(x) = y$. 
\begin{thm}  \label{thm-bing}
Let $X$ be    a homogeneous, circle-like
continuum that contains an arc. Then either $X$ is homeomorphic to the circle $\mS^1$, or to    an inverse limit of coverings of   $\mS^1$.
\end{thm}
This results inspired the subsequent works  by McCord \cite{McCord1965}, Thomas \cite{EThomas1973}, Hagopian \cite{Hagopian1977}, Mislove and Rogers \cite{MR1989}, and Aarts, Hagopian and Oversteegen \cite{AHO1991}, all for $1$-dimensional flow spaces. 

In this work, we give extensions of Theorem~\ref{thm-bing} to continua with higher dimensional arc-components. We first recall two notions which are required to formulate our results.  A weak solenoid $\cS_\cP$ is the inverse limit space of a sequence of   covering maps of finite degree greater than one,
\begin{equation}\label{eq-presentation}
\cP = \{p_{\ell +1} \colon M_{\ell +1} \to M_{\ell} \mid \ell \geq 0\},
\end{equation}
   where $M_\ell$ is a compact connected manifold without boundary. The collection of maps $\cP$ is called a \emph{presentation} for $\cS_\cP$.  
  A weak solenoid $\cS_\cP$ is \emph{regular} if the presentation $\cP$ can be chosen so that for each $\ell \geq 0$, the composition $p_{\ell}^0 \equiv p_1 \circ \cdots p_{\ell} \colon M_\ell \to M_0$ is a regular covering map; that is, the fundamental group of $M_\ell$
injects onto a normal subgroup of the fundamental group of $M_0$ under the
map induced by the covering projection $p_{\ell}^0$.   A weak solenoid which is not regular is said to be \emph{irregular}.   A \emph{Vietoris solenoid} \cite{vanDantzig1930,Vietoris1927}  is   a $1$-dimensional  regular solenoid, where each $M_{\ell}$ is a circle, as arises in the conclusion of Theorem~\ref{thm-bing}.

A matchbox manifold    $\fM$ is a continuum  equipped with a decomposition $\F$ into leaves of constant dimension, 
 so that the pair $(\fM, \F)$ is a foliated space  in the sense of \cite{MS2006}, for which the local transversals to the foliation  are totally disconnected. In particular,    
the leaves of $\F$ are the path connected components of $\fM$.   A matchbox manifold with $2$-dimensional leaves is a lamination by surfaces in the sense of Ghys  \cite{Ghys1999} and Lyubich and Minsky \cite{LM1997}, while Sullivan called them ``solenoidal spaces'' in \cite{Sullivan2014,Verjovsky2014}. The terminology ``matchbox manifold'' follows  the usage introduced in  \cite{AM1988,AO1991,AO1995}. 
    A Vietoris solenoid is a $1$-dimensional matchbox manifold, and more generally, McCord showed in \cite{McCord1965} that 
  $n$-dimensional solenoids   are examples of $n$-dimensional  matchbox manifolds.
 
 Next, recall a result of the first two authors. A matchbox manifold is said to be \emph{equicontinuous} if the holonomy pseudogroup $\cGF$ associated to the foliation $\F$ (see Section~\ref{sec-holonomy}) defines an equicontinuous  action on its transversal space, as defined in Definition~\ref{def-equicontinuous}.
 \begin{thm} \cite[Theorems~1.2 and 1.4]{ClarkHurder2013} \label{thm-CH2013}
Let $\fM$ be an equicontinuous   matchbox manifold. Then $\fM$ is
homeomorphic to a  weak solenoid, and in particular is manifold-like.  Moreover, if $\fM$ is homogeneous, then $\fM$ is homeomorphic to a regular solenoid.
\end{thm}
 
 The main result  of  this paper, as follows, yields a generalization of   Theorem~\ref{thm-bing} to higher dimensional matchbox manifolds. 
\begin{thm}\label{thm-main1}
A manifold-like matchbox manifold  $\fM$ is equicontinuous.
\end{thm}

  Theorems~\ref{thm-CH2013} and \ref{thm-main1} then yield the following  partial    converse to Theorem~\ref{thm-CH2013}:

\begin{cor}\label{cor-weak}
A    manifold-like matchbox manifold $\fM$ is homeomorphic to a  weak  solenoid. 
\end{cor}

The hypothesis that a matchbox manifold $\fM$ is manifold-like does not imply that $\fM$ is homogeneous, as the ``discriminant obstruction'' to homogeneity for a weak solenoid is supported in arbitrarily small open neighborhoods of points in $\fM$ . The discriminant invariant was introduced and studied in the works in \cite{DHL2016a,DHL2016b,DHL2016c,HL2017a}.
  
The remainder of this paper is organized as follows. In Section~\ref{sec-concepts} we recall the definitions of foliated spaces and matchbox manifolds, and give some of their basic properties. Particular care is taken to introduce various metric estimates related to the geometry of the leaves of the foliation, and to its dynamical properties. 
In Section~\ref{sec-holonomy} we recall the construction of the holonomy along leafwise paths. 

In Section~\ref{sec-pathlifting}, we prove Theorem~\ref{thm-main1}, using the path lifting property for $\e$-maps from a matchbox manifold to a compact manifold. This is the key technical tool, which is used to show that the foliation $\F$ on $\fM$ must be equicontinuous.   

The Appendix~\ref{appendix-A} contains the proof of a technical result, Proposition~\ref{prop-coverage} below. The proof of this result does not appear to be in the literature, and may even have a simpler proof than is given here. However, the result is essential for the proof of  Theorem~\ref{thm-main1}, so is included for completeness.
 
\eject

\section{Foliated spaces and matchbox manifolds} \label{sec-concepts}

We recall some   background concepts used in the proof of our main theorems.  

\subsection{Matchbox manifolds}\label{subsec-fs} 
We first recall the definition  of a  matchbox manifold.

\begin{defn} \label{def-fs}
A \emph{matchbox manifold of dimension $n$} is a   continuum $\fM$, such that  there exists a compact, separable, totally disconnected metric space $\fX$, and
for each $x \in \fM$ there is a compact subset $\fT_x \subset \fX$, an open subset $U_x \subset \fM$, and a homeomorphism $\vp_x \colon \oU_x \to [-1,1]^n \times \fT_x$ defined on the closure $\oU_x$ in $\fM$,
  such that $\vp_x(x) = (0, w_x)$ where $w_x \in int(\fT_x)$. Moreover, it is assumed that each $\vp_x$  admits an extension to a foliated homeomorphism
$\whvarp_x \colon \whU_x \to (-2,2)^n \times \fT_x$ where $\whU_x \subset \fM$ is an open subset such that $\oU_x \subset \whU_x$.
The space $\fT_x$   is  called the \emph{local transverse model} at $x$.  
\end{defn}

The assumption that the transversals $\fT_x$ are totally disconnected implies that the local charts $\vp_x$ satisfy the compatibility axioms of foliation charts for a foliated space, as   in  \cite{CandelConlon2000,ClarkHurder2013,MS2006}.

Let $\pi_x \colon \oU_x \to \fT_x$ denote the composition of $\vp_x$ with projection onto the second factor.

Also introduce the transversal maps $\tau_x \colon \fT_x \to \cT_x \subset \fM$, defined for $w \in \fT_x$ by $\tau_x(w) = \vp_x^{-1}(0 , w)$.
The subspace $\cT_x$ is   given the   metric $d_{\cT_x}$ which is the restriction of  the metric $\dM$.

For $w \in \fT_x$ the set $\cP_x(w) = \pi_x^{-1}(w) \subset \oU_x$ is called a \emph{plaque} for the coordinate chart $\vp_x$. We adopt the notation, for $z \in \oU_x$, that $\cP_x(z) = \cP_x(\pi_x(z))$, so that $z \in \cP_x(z)$. Note that each plaque $\cP_x(w)$ is given the topology so that the restriction $\vp_x \colon \cP_x(w) \to [-1,1]^n \times \{w\}$ is a homeomorphism. Then $int (\cP_x(w)) = \vp_x^{-1}((-1,1)^n \times \{w\})$.

Let $U_x = int (\oU_x) = \vp_x^{-1}((-1,1)^n \times int(\fT_x))$.
Note that if $z \in U_x \cap U_y$, then $int(\cP_x(z)) \cap int( \cP_y(z))$ is an open subset of both
$\cP_x(z) $ and $\cP_y(z)$.
The collection of sets
$$\cV = \{ \vp_x^{-1}(V \times \{w\}) \mid x \in \fM ~, ~ w \in \fT_x ~, ~ V \subset (-1,1)^n ~ {\rm open}\}$$
forms the basis for the \emph{fine topology} of $\fM$. The connected components of the fine topology are called leaves, and define the foliation $\F$ of $\fM$.
 In particular,    the leaves of the foliation $\F$ of $\fM$  are the path-connected components of $\fM$.
For $x \in \fM$, let $L_x \subset \fM$ denote the leaf of $\F$ containing $x$.

\begin{defn} \label{def-sfs}
A \emph{smooth matchbox manifold} is a   space $\fM$ as above, such that there exists a choice of local charts $\vp_x \colon \oU_x \to [-1,1]^n \times \fT_x$ such that for all $x,y \in \fM$ with $z \in U_x \cap U_y$, there exists an open set $z \in V_z \subset U_x \cap U_y$ such that $\cP_x(z) \cap V_z$ and $\cP_y(z) \cap V_z$ are connected open sets, and the composition
$$\psi_{x,y;z} \equiv \vp_y \circ \vp_x ^{-1}\colon \vp_x(\cP_x (z) \cap V_z) \to \vp_y(\cP_y (z) \cap V_z)$$
is a smooth map, where $\vp_x(\cP_x (z) \cap V_z) \subset \mR^n \times \{w\} \cong \mR^n$ and $\vp_y(\cP_y (z) \cap V_z) \subset \mR^n \times \{w'\} \cong \mR^n$. The leafwise transition maps $\psi_{x,y;z}$ are assumed to depend continuously on $z$ in the $C^{\infty}$-topology.
\end{defn}

A map $f \colon \fM \to \mR$ is said to be \emph{smooth} if for each flow box
$\vp_x \colon \oU_x \to [-1,1]^n \times \fT_x$ and $w \in \fT_x$ the composition
$y \mapsto f \circ \vp_x^{-1}(y, w)$ is a smooth function of $y \in (-1,1)^n$, and depends continuously on $w$ in the $C^{\infty}$-topology on maps of the plaque coordinates $y$. As noted in \cite{MS2006} and \cite[Chapter 11]{CandelConlon2000}, this allows one to define smooth partitions of unity, vector bundles, and tensors for smooth foliated spaces. In particular, one can define leafwise Riemannian metrics. We recall a standard result, whose proof for foliated spaces can be found in \cite[Theorem~11.4.3]{CandelConlon2000}.
\begin{thm}\label{thm-riemannian}
Let $\fM$ be a smooth matchbox manifold. Then there exists a leafwise Riemannian metric for $\F$, such that for each $x \in \fM$, the leaf $L_x$ inherits the structure of a complete Riemannian manifold with bounded geometry, and the Riemannian metric and its covariant derivatives depend continuously on $x$ .  
\end{thm}
Bounded geometry implies, for example, that for each $x \in \fM$, there is a leafwise exponential map
$\exp^{\F}_x \colon T_x\F \to L_x$ which is a surjection, and the composition $\exp^{\F}_x \colon T_x\F \to L_x \subset \fM$ depends continuously on $x$ in the compact-open topology on maps.
 All matchbox manifolds  are assumed to be smooth with a given leafwise Riemannian metric, and with a fixed choice of    metric $\dM$ on $\fM$.

\subsection{Metric estimates}\label{subsec-metric}
We formulate some relations between the metric properties of a matchbox manifold $\fM$ and the metric properties of the leaves of $\F$. These technical conditions are used in studying the dynamics and geometry of these spaces.

For $x \in \fM$ and $\e > 0$, let $D_{\fM}(x, \e) = \{ y \in \fM \mid \dM(x, y) \leq \e\}$ be the closed $\e$-ball about $x$ in $\fM$, and $B_{\fM}(x, \e) = \{ y \in \fM \mid \dM(x, y) < \e\}$ the open $\e$-ball about $x$.

Recall that $\dX$ denotes the metric on the   space $\fX$ in Definition~\ref{def-mm}.  For $w \in \fX$ and $\e > 0$, let $D_{\fX}(w, \e) = \{ w' \in \fX \mid \dX(w, w') \leq \e\}$ be the closed $\e$-ball about $w$ in $\fX$, and let $B_{\fX}(w, \e) = \{ w' \in \fX \mid \dX(w, w') < \e\}$ be the open $\e$-ball about $w$.

Each leaf $L \subset \fM$ has a complete path-length metric, induced from the leafwise Riemannian metric:
$$\dF(x,y) = \inf \left\{\| \gamma\| \mid \gamma \colon [0,1] \to L ~{\rm is ~ piecewise ~~ C^1}~, ~ \gamma(0) = x ~, ~ \gamma(1) = y ~, ~ \gamma(t) \in L \quad \forall ~ 0 \leq t \leq 1\right\}$$
  where $\| \gamma \|$ denotes the path-length of the piecewise $C^1$-curve $\gamma(t)$. If $x,y \in \fM$   are not on the same leaf, then set $\dF(x,y) = \infty$.

  For each $x \in \fM$ and $r > 0$, let $D_{\F}(x, r) = \{y \in L_x \mid \dF(x,y) \leq r\}$.

For each $x \in \fM$, the  {Gauss Lemma} implies that there exists $\lambda_x > 0$ such that $D_{\F}(x, \lambda_x)$ is a \emph{strongly convex} subset for the metric $\dF$. That is, for any pair of points $y,y' \in D_{\F}(x, \lambda_x)$ there is a unique shortest geodesic segment in $L_x$ joining $y$ and $y'$ and  contained in $D_{\F}(x, \lambda_x)$. This standard concept of Riemannian geometry is discussed in detail in  \cite{BC1964}, and in \cite[Chapter 3, Proposition 4.2]{doCarmo1992}. 
Then for all $0 < \lambda < \lambda_x$ the disk $D_{\F}(x, \lambda)$ is also strongly convex. 
The leafwise metrics for $\F$  constructed in the proof of Theorem~\ref{thm-riemannian} have uniformly bounded geometry, and the  first and second order covariant derivatives of the metrics depend continuously on the point $x \in \fM$, so by the compactness of $\fM$, we obtain:
\begin{lemma}\label{lem-stronglyconvex}
There exists $\lF > 0$ such that for all $x \in \fM$, $D_{\F}(x, \lF)$ is strongly convex.
\end{lemma}

\subsection{Regular coverings} \label{subsec-rc}
We next formulate  the definition of a \emph{regular covering} of a matchbox manifold $\fM$.
It follows from standard considerations (see \cite{ClarkHurder2013}) that a matchbox manifold admits a covering by foliation charts which satisfies additional regularity conditions.

\begin{prop}\label{prop-regular}\cite{ClarkHurder2013}
For a smooth foliated space $\fM$, given $\eM > 0$, there exist   $\lF>0$ and a choice of local charts $\vp_x \colon \oU_x \to [-1,1]^n \times \fT_x$ with the following properties:  For each $x \in \fM$,
\begin{enumerate}
\item  $U_x \equiv int(\oU_x) = \vp_x^{-1}\left( (-1,1)^n \times \fT_x \right)$, with   $\oU_x \subset B_{\fM}(x, \eM)$.
\item  The plaques of $\vp_x$ are   strongly convex for the metric $\dF$ with diameter less than $\lF$. 
\end{enumerate}
\end{prop}

By a standard argument, there exists a finite collection  $\{x_1, \ldots , x_{\nu}\} \subset \fM$ where $\vp_{x_i}(x_i) = (0 , w_{x_i})$ for $w_{x_i} \in \fX$,  and regular foliation charts $\vp_{x_i} \colon \oU_{x_i} \to [-1,1]^n \times \fT_{x_i}$
satisfying the conditions of Proposition~\ref{prop-regular},   which form an open covering of $\fM$. Relabel the various maps and spaces accordingly, so that
$\oU_{i} = \oU_{x_i}$ and $\vp_{i} = \vp_{x_i}$. Accordingly, label the transverse spaces $\fT_i = \fT_{x_i}$ and the projection maps $\pi_i = \pi_{x_i} \colon \oU_i \to \fX_i$.
  Then the image $\pi_i(U_i \cap U_j) = \fT_{i,j} \subset \fT_i$ is a clopen subset for all $1 \leq i, j \leq \nu$.   

We also then have the transversal mappings 
$\tau_i \colon \fT_i \to \cT_i \subset \fM$ for $1 \leq i \leq \nu$.

A \emph{regular covering} of $\fM$ is a finite covering $\cU = \{U_{1}, \ldots , U_{\nu}\}$ such that for each $1 \leq i \leq \nu$ there is a foliated coordinate map 
$\ds \vp_i \colon  U_i \to (-1,1)^n \times \fT_i$  which satisfies  the   regularity  conditions in Proposition~\ref{prop-regular}.
  We assume in the following that a   regular foliated covering  of $\fM$   has been chosen.

\subsection{More metric estimates}\label{subsec-moremetric}
 We introduce lower and upper bounds on the diameters of balls in the leaves of $\F$ with respect to the ambient metric $\dM$ on $\fM$. To assist with the notation,  we use the convention that $\lambda > 0$ will  denote a small leafwise distance, and $\e$ will denote a small distance in $\fM$. Later when we introduce the target manifold $M$, we let  $\delta$   denote a small distance in $M$. 
 
 For $x \in \fM$ and $\e > 0$,   let $D_{\F}(\dM, x, \e) \subset L_x$ denote the connected component containing $x$ of the intersection $L_x \cap D_{\fM}(x, \e)$.  
Define the   continuous functions
\begin{eqnarray}
\rho(\dF,\dM, x,\e) & = & \max \ \{\dF(x',x) \mid x' \in D_{\F}(\dM, x, \e)\} \label{eq-rhoFM}\\
\rho(\dM,\dF, x,\lambda) & = & \max \ \{\dM(x',x) \mid x' \in D_{\F}(x, \lambda)\} \ . \label{eq-rhoMF}
\end{eqnarray}
Then for all $x \in \fM$,  $\e > 0$ and $\lambda > 0$,  we have
\begin{equation}
D_{\F}(\dM, x, \e) \subset D_{\F}(x, \rho(\dF,\dM, x,\e)) ~ {\rm and}~  D_{\F}(x, \lambda) \subset D_{\fM}(x, \rho(\dM,\dF, x,\lambda)) \ .
\end{equation}
As $\fM$ is compact, we can then define   the  increasing   functions of $\e > 0$ and $\lambda > 0$,
 \begin{eqnarray}
\rho(\dF,\dM, \e) & = & \max \ \{\rho(\dF,\dM, x,\e) \mid x \in \fM\} \\
\rho(\dM,\dF, \lambda) & = & \max  \ \{\rho(\dM,\dF, x,\lambda) \mid x \in \fM\} \ .
\end{eqnarray}
Moreover, for $0 < \lambda' < \lambda \leq \lF$ we have the strict inclusion $B_{\F}(x,\lambda') \subset B_{\F}(x,\lambda)$, and thus  
the function $\rho(\dM,\dF, \lambda)$ is strictly increasing for $0 < \lambda <   \lF$.

Let $\e_{\F}^* = \max \ \{ \e \mid \rho(\dF, \dM, \e) \leq \lF \}$, so that $D_{\F}(\dM, x, \e_{\F}^*) \subset D_{\F}(x, \lF)$ for all $x \in \fM$.

Introduce  the continuous function  
$\lambda_{\F}(\e)$   which is the inverse of $\rho(\dM,\dF, \lambda)$, for  $0 < \lambda \leq \lF$. Thus, $\lambda_{\F}(\e)$ is the largest radius $\lambda \leq \lF$ such that the   disk  $D_{\F}(x, \lambda)$ in the leafwise metric is contained in the ball $B_{\fM}(x,\e)$ for all $x \in \fM$. Combining the above definitions, we obtain that    for   all $x \in \fM$, 
\begin{equation}\label{eq-metriccont}
  D_{\F}(x, \lambda_{\F}(\e))   \subset  D_{\F}(\dM, x, \e) \subset  D_{\F}(x, \lF) \cap  D_{\fM}(x, \e)  \ .
\end{equation}
  
Choose $\e_0 > 0 $ so that  $\rho(\dF,\dM, \e_0) \leq \lF/2$,  and set $\lambda_0 = \lambda_{\F}(\e_0)$.
Then by the definition of $\rho(\dF,\dM, \e_0)$ and the inclusions         \eqref{eq-metriccont},   for all $x \in \fM$ we have the inclusions
\begin{equation}
D_{\F}(x, \lambda_0 )   \subset     D_{\F}(\dM, x, \e_0) \subset  D_{\F}(x, \lF/2)  \cap    B_{\fM}(x, \e_0) \ .
\end{equation}

Next, choose $\e_1 > 0$ so that  $\rho( \dF, \dM, \e_1) \leq \lambda_0/10$, and let $\lambda_1 = \lF(\e_1)$.
Then  for all $x \in \fM$,  
\begin{equation}\label{eq-disklambda0}
D_{\F}(x, \lambda_1)   \subset  D_{\F}(\dM, x, \e_1)  \subset D_{\F}(x, \lambda_0/10)  \cap    B_{\fM}(x, \e_1) \ .
\end{equation}
 This choice of $\e_1$ will be recalled in Section~\ref{sec-pathlifting} and Appendix~\ref{appendix-A}.

   A matchbox manifold $\fM$ is \emph{minimal} if every leaf of $\F$ is dense.

\section{Holonomy} \label{sec-holonomy}

The holonomy pseudogroup of a smooth foliated manifold $(M, \F)$ generalizes the induced dynamical systems  associated to a section of a flow. The holonomy pseudogroup for a matchbox manifold $(\fM, \F)$ is defined analogously to the smooth case.

 \subsection{The foliation pseudo-star group}\label{subsec-psg}
  Let $\cU = \{ \vp_{i}  \colon \oU_i \to [-1,1]^n \times \fT_i  \mid 1 \leq i \leq \nu\}$ be a regular covering of   $\fM$ as in Section~\ref{subsec-rc}. 
   A pair of indices $(i,j)$, $1 \leq i,j \leq \nu$, is said to be \emph{admissible} if the \emph{open} coordinate charts satisfy $U_i \cap U_j \ne \emptyset$.
For $(i,j)$ admissible, define clopen subsets $\fD_{i,j} = \pi_i(U_i \cap U_j) \subset \fT_i \subset \fX$.  The convexity of foliation charts imply that plaques are either disjoint, or have connected intersection. This implies that there is a well-defined homeomorphism $h_{j,i} \colon \fD_{i,j} \to \fD_{j,i}$ with domain $\fD(h_{j,i}) = \fD_{i,j}$ and range $R(h_{j,i}) = \fD_{j,i}$.

The maps $\cGF^{(1)} = \{h_{j,i} \mid (i,j) ~{\rm admissible}\}$ are the transverse change of coordinates defined by the foliation charts. By definition they satisfy $h_{i,i} = Id$, $h_{i,j}^{-1} = h_{j,i}$, and if $U_i \cap U_j\cap U_k \ne \emptyset$ then $h_{k,j} \circ h_{j,i} = h_{k,i}$ on their common domain of definition. The \emph{holonomy pseudogroup} $\cGF$ of $\F$ is the topological pseudogroup modeled on $\fX$ generated by   the elements of $\cGF^{(1)}$. The elements of $\cGF$ have a standard description in terms of the ``holonomy along paths'', which we next describe.

A sequence $\cI = (i_0, i_1, \ldots , i_{\alpha})$ is \emph{admissible}, if each pair $(i_{\ell -1}, i_{\ell})$ is admissible for $1 \leq \ell \leq \alpha$, and the composition
\begin{equation}\label{eq-defholo}
h_{\cI} = h_{i_{\alpha}, i_{\alpha-1}} \circ \cdots \circ h_{i_1, i_0}
\end{equation}
has non-empty domain.
The domain $\fD_{\cI}$ of $h_{\cI}$ is the \emph{maximal clopen subset} of $\fD_{i_0} \subset \fT_{i_0}$ for which the compositions are defined.

 For the study of the dynamical properties of $\F$, it is necessary  to introduce   the collection of maps $\cGF^* \subset \cGF$, defined as follows.
Given any open subset $U \subset \fD_{\cI}$ we obtain a new element $h_{\cI} | U \in \cGF$ by restriction. Then set
\begin{equation}\label{eq-restrictedgroupoid}
\cGF^* = \left\{ h_{\cI} | U \mid \cI ~ {\rm admissible} ~ \& ~ U \subset \fD_{\cI} \right\} \subset \cGF ~ .
\end{equation}
  That is,  $\cGF^*$   consists of all possible restrictions of homeomorphisms of the form \eqref{eq-defholo} to open subsets of their domains.  
 However,   in the definition of    $\cGF^*$   one does not allow arbitrary unions of local homeomorphisms, unless such homeomorphisms can be obtained by restrictions to open subsets of maximal domains of words in the elements in $\cG_0$. 
  The collection of maps $\cGF^*$    is closed under the operations of compositions, taking inverses, and restrictions to open sets, and  is called a \emph{$\psg$} in the literature \cite{Matsumoto2010}.  
  
 For $g \in \cGF^*$ denote its domain by $\fD(g) \subset \fX$, then its  range   is the clopen set $\fR(g) = g(\fD(g))  \subset \fX$.

 \subsection{Admissible chains}\label{subsec-admissible}

Given an admissible sequence $\cI = (i_0, i_1, \ldots , i_{\alpha})$ and any $0 \leq \ell \leq \alpha$,     the truncated sequence $\cI_{\ell} = (i_0, i_1, \ldots , i_{\ell})$ is again admissible, and we introduce the holonomy map defined by the composition of the first $\ell$ generators appearing in $h_{\cI}$,
\begin{equation}\label{eq-pcmaps}
h_{\cI_{\ell}} = h_{i_{\ell} , i_{\ell -1}} \circ \cdots \circ h_{i_{1} , i_{0}}~.
\end{equation}
Given $w \in \fD(h_{\cI})$ we adopt the notation $w_{\ell} = h_{\cI_{\ell}}(w) \in \fT_{i_{\ell}}$. So $w_0 = w$ and
$h_{\cI}(w) = w_{\alpha}$.

Given $w \in \fD(h_{\cI})$, let $x_0 = \tau_{i_0}(w_0) \in L_{x_0}$. Introduce the \emph{plaque chain}
\begin{equation}\label{eq-plaquechain}
\cP_{\cI}(w) = \{\cP_{i_0}(w_0), \cP_{i_1}(w_1), \ldots , \cP_{i_{\alpha}}(w_{\alpha}) \} ~ .
\end{equation}

 For each $1 \leq i \leq \nu$, define    $\cT_{i} =  \vp_i^{-1}(0 , \fT_i) \subset    \oU_i \subset \fM$ which is a compact local transversal to $\F$.

Intuitively, a plaque chain $\cP_{\cI}(w)$ is a sequence of successively overlapping convex ``tiles'' in $L_0$ starting at $x_0 = \tau_{i_0}(w_0)$, ending at
$x_{\alpha} = \tau_{i_{\alpha}}(w_{\alpha})$, and with each $\cP_{i_{\ell}}(w_{\ell})$ ``centered'' on the point $x_{\ell} = \tau_{i_{\ell}}(w_{\ell})$.
Recall that $\cP_{i_{\ell}}(x_{\ell}) = \cP_{i_{\ell}}(w_{\ell})$, so we also adopt the notation $\cP_{\cI}(x) \equiv \cP_{\cI}(w)$.

A \emph{leafwise path}  is a continuous map $\gamma \colon [0,1] \to \fM$ such that there is a leaf $L$ of $\F$ for which $\gamma(t) \in L$ for all $0 \leq t \leq 1$.
  In the following, we will assume that all paths are  piecewise differentiable.
  
 Let $\gamma$ be a leafwise path, and $\cI$ be an admissible sequence. For $w \in \fD(h_{\cI})$, we say that $(\cI , w)$ \emph{covers} $\gamma$,
if the domain of $\gamma$ admits   a partition $0 = s_0 < s_1 < \cdots < s_{\alpha} = 1$ such that  the plaque chain $\cP_{\cI}(w_0) = \{\cP_{i_0}(w_0), \cP_{i_1}(w_1), \ldots , \cP_{i_{\alpha}}(w_{\alpha}) \}$ satisfies
\begin{equation}\label{eq-cover}
\gamma([s_{\ell} , s_{\ell + 1}]) \subset int (\cP_{i_{\ell}}(w_{\ell}) )~ , ~ 0 \leq \ell < \alpha, \quad \& \quad \gamma(1) \in int( \cP_{i_{\alpha}}(w_{\alpha})).
\end{equation}
  The map $h_{\cI}$ is said to define the holonomy of $\F$ along the path $\gamma$,  and satisfies $h_{\cI}(w_0) = \pi_{i_{\alpha}}(\gamma(1))$.

Given two admissible sequences, $\cI = (i_0, i_1, \ldots, i_{\alpha})$ and
$\cJ = (j_0, j_1, \ldots, j_{\beta})$, such that both $(\cI, w_0)$ and $(\cJ, v_0)$ cover the leafwise path $\gamma \colon [0,1] \to \fM$, then
$$\gamma(0) \in int( \cP_{i_0}(w_0)) \cap int( \cP_{j_0}(v_0)) \quad , \quad \gamma(1) \in int(\cP_{i_{\alpha}}(w_{\alpha})) \cap int( \cP_{j_{\beta}}(v_{\beta}) )$$ Thus both $(i_0 , j_0)$ and $(i_{\alpha} , j_{\beta})$ are admissible, and
$v_0 = h_{j_{0} , i_{0}}(w_0)$, $w_{\alpha} = h_{i_{\alpha} , j_{\beta}}(v_{\beta})$.

The proof of the following standard observation can be found in \cite{ClarkHurder2013}.
\begin{prop}\label{prop-copc}\cite{ClarkHurder2013}
The maps $h_{\cI}$ and
$\ds h_{i_{\alpha} , j_{\beta}} \circ h_{\cJ} \circ h_{j_{0} , i_{0}}$
agree on their common domains.
\end{prop}

Two leafwise paths $\gamma , \gamma' \colon [0,1] \to \fM$ are \emph{homotopic} if there exists a family of leafwise paths $\gamma_s \colon [0,1] \to \fM$ with $\gamma_0 = \gamma$ and $\gamma_1 = \gamma'$. We are most interested in the special case when $\gamma(0) = \gamma'(0) = x$ and $\gamma(1) = \gamma'(1) = y$. Then $\gamma$ and $\gamma'$ are \emph{homotopic relative endpoints}, or \emph{endpoint-homotopic}, 
if they are homotopic with $\gamma_s(0) = x$ for all $0 \leq s \leq 1$, and similarly
$\gamma_s(1) = y$ for all $0 \leq s \leq 1$. Thus, the family of curves $\{ \gamma_s(t) \mid 0 \leq s \leq 1\}$ are all contained in a common leaf $L_{x}$. We then have the following result.

\begin{lemma}\cite{ClarkHurder2013} \label{lem-homotopic}
Let $\gamma, \gamma' \colon [0,1] \to \fM$ be endpoint-homotopic leafwise paths. Then their holonomy maps $h_{\gamma}$ and $h_{\gamma'}$ agree on some open subset $U \subset \fD(h_{\gamma}) \cap \fD(h_{\gamma'}) \subset \fT_*$. 
\end{lemma}

Finally, we   recall the     definition of an  {equicontinuous} pseudogroup.
\begin{defn} \label{def-equicontinuous}
The action of the    pseudogroup $\cGF$ on $\fX$ is  \emph{equicontinuous} if for all $\epsilon > 0$, there exists $\delta > 0$ such that for all $g \in \cGF^*$, if $w, w' \in \fD(g)$ and $\dX(w,w') < \delta$, then $\dX(g(w), g(w')) < \epsilon$.
Thus, $\cGF^*$ is equicontinuous as a family of local group actions.
\end{defn}
 Further   properties of the pseudogroup $\cGF$ for a matchbox manifold are discussed in  \cite{ClarkHurder2013,CHL2014,Hurder2014}.

\section{Equicontinuous holonomy} \label{sec-pathlifting}
 
In this section, we give the proof of Theorem~\ref{thm-main1}. A key point in the proof  is based on   the path lifting property for $\e$-maps between  a matchbox manifold $\fM$  and a target manifold $M$. 
The philosophy of path lifting is   folklore, as observed by Bob Edwards in his 1978  ICM address on $\e$-approximations and homeomorphisms \cite[Section~4]{Edwards1980}. 
We develop  this technique in the context of  matchbox manifolds, making use of   standards results for $\e$-maps along  with   properties of the holonomy maps.
  
 Let $\fM$ be a   manifold-like matchbox manifold.     Let $\cU = \{ \vp_{i}  \colon \oU_i \to [-1,1]^n \times \fT_i  \mid 1 \leq i \leq \nu\}$ be a regular covering of   $\fM$ as in Section~\ref{subsec-rc}.
 Let $\cGF^*$ be the $\psg$ associated to the regular covering $\cU$ as in    Section~\ref{subsec-psg}.
We must show that the conditions of Definition~\ref{def-equicontinuous} are satisfied: 
  given $\epsilon > 0$, we must show there exists $\delta > 0$ such that  for each admissible chain $\cI$ with holonomy $h_{\cI}$, if $w, w' \in \fD(h_{\cI}) \subset \fX$ and $\dX(w,w') < \delta$, then $\dX(h_{\cI}(w), h_{\cI}(w')) < \e$.

Recall that $\eU > 0$ denotes a Lebesgue number for the covering $\cU$ of $\fM$. That is, for each $x \in \fM$, there exists some index $1 \leq i \leq \nu$ such that $B_{\fM}(x, \eU) \subset U_{i}$.

Since  $\oU_i$   is compact for each $1 \leq i \leq \nu$, there exists a   uniform modulus of continuity function   $\rho_{\cU}(\e) > 0$ for the projections $\pi_i$:  let $\e > 0$, then   $\rho_{\cU}(\e)$ is the largest value  such that 
\begin{equation}\label{eq-unifcontU}
B_{\fM}(x,\rho_{\cU}(\e)) \cap \oU_i \subset \pi_i^{-1}(B_{\fX_i}(\pi_i(x), \e)) ~ {\rm for ~ all}~  x \in  \oU_i \ . 
\end{equation}
For $\e_1$ as defined in Section~\ref{subsec-metric}, define $\eF$ by 
\begin{equation}\label{eq-eF}
0 < \eF = \min \ \{\rho_{\cU}(\e), \e_1/2 ,  \eU/4  \} \ .
\end{equation}

 \subsection{Continuity estimates}\label{subsec-contest}

Choose an $\eF$-map $f \colon \fM \to M$ onto the compact topological manifold $M$, where for simplicity we omit the subscript $\eF$ in the notation for $f$.
Let $d_M$ be a metric on $M$. 
For   $w \in M$ and $\delta > 0$, let   $B_{M}(w, \delta) =\{w' \in M \mid d_M(w',w) < \delta\}$  denote the open disk in $M$ of radius $\delta$, and     $D_{M}(w, \delta) =\{w' \in M \mid d_M(w',w) \leq \delta\}$   denote the closed disk in $M$ of radius $\delta$.

 Assume that $d_M$ is chosen so that there exists 
a constant    $\delta_M >0$  such that for all $w \in M$ and $0 < \delta \leq \delta_M$, the   disk $D_{M}(w, \delta)$ is homeomorphic to a disk in $\mR^n$. For example, if $M$ is a smooth Riemannian manifold, then   let $\delta_M > 0$ be   such that each   disk  $D_{M}(w, \delta_M)$ is strongly convex.

Since $\fM$ is compact, there exists a   uniform modulus of continuity function   $\e_f(\delta) > 0$ for $f$:  for $\delta > 0$, the constant $\e_f(\delta)$ is the largest value  such that 
\begin{equation}\label{eq-unifcont}
B_{\fM}(x,\e_f(\delta)) \subset f^{-1}(B_{M}(f(x), \delta)) ~ {\rm for ~ all}~  x \in \fM\ . 
\end{equation}
Let     $\delta_f^* > 0$ be the largest radius such that $\e_f(\delta) \leq \e_{\F}$ for all $0 < \delta \leq \delta_f^*$.

  Set $\lambda_{\F,f}(\delta) = \lambda_{\F}(\e_f(\delta))$, which is
well-defined for $0 < \delta \leq \delta_f^*$. Recall from Section
\ref{subsec-moremetric} that $\lambda_{\F,f}(\delta)$ is then Ê  the
largest radius $\lambda \leq \lF$ such that the ÊÊdisk Ê$D_{\F}(x,
\lambda)$ in the leafwise metric is contained in the ball
$B_{\fM}(x,\e_f(\delta))$ for all $x \in \fM$.

 Combining \eqref{eq-metriccont} and \eqref{eq-unifcont} we obtain a leafwise modulus of continuity for $f$:  
 \begin{equation}\label{eq-unifleafcont}
  D_{\F}(x, \lambda_{\F,f}(\delta) )   \subset  B_{\fM}(x,\e_f(\delta)) \subset f^{-1}(B_{M}(f(x), \delta)) \quad {\rm for ~ all}~ 0 < \delta \leq \delta_f^* ~, ~ x \in \fM\ .   
\end{equation}
 As $f$ is an $\eF$-map, we have that  $f^{-1}(f(x)) \subset B_{\fM}(x,\eF)$ for all $x \in \fM$. 
As $M$ is compact and $f^{-1}(f(x))$ is a compact set with diameter at most $\eF$ for all $x \in \fM$,  there exists  $0 < \delta_1   \leq \delta_M/10$ so  that  for all $x \in \fM$, we have 
 \begin{equation}\label{eq-diskdelta1}
 f^{-1}(D_M(f(x),\delta_1)) \subset B_{\fM}(x,2\eF) = B_{\fM}(x,\e_1) \ .
 \end{equation}

Let $\lambda_2 =  \lambda_{\F,f}(\delta_1)$; that is, $\lambda_2$ is then Êis the largest radius $\lambda \leq \lF$
such that the ÊÊdisk Ê$D_{\F}(x, \lambda)$ in the leafwise metric is
contained in the ball $B_{\fM}(x,\delta_1)$ for all $x \in \fM$.

Then for all $x \in \fM$,  by \eqref{eq-unifleafcont},  \eqref{eq-disklambda0},  and the choice of $\delta_1$ we have
  \begin{equation}\label{eq-diameter}
D_{\F}(x, \lambda_2 )   \subset f^{-1}(D_M(f(x),\delta_1)) \cap D_{\F}(x,\lF)   \subset D_{\fM}(x,\e_1) \cap D_{\F}(x,\lF)   \subset D_{\F}(x,\lambda_0/10) \ .
\end{equation}

 Finally, we     require a basic result  concerning $\e$-maps on matchbox manifolds, which is a type of converse to the inclusions in \eqref{eq-diameter}, and whose   proof  is   in the spirit of the work by Eilenberg  \cite{Eilenberg1938}. The proof of the following  is deferred to  Appendix~\ref{appendix-A}.

\begin{prop}\label{prop-coverage}
Let $\fM$ be a matchbox manifold with leafwise Riemannian metric on $\F$. Then 
there exists $\eF > 0$  such that,    if $f \colon \fM \to M$ is an $\eF$-map to a compact manifold $M$, then for      $x_0 \in \fM$ with $w_0 = f(x_0)$,  we have  $B_M(w_0,\delta_1) \subset f(D_{\F}(x_0, \lF/2 ))$.
\end{prop}

   \subsection{Local lifting property}\label{subsec-extensions}

  We next establish  a   technical result used in the proof of Theorem~\ref{thm-main1}. 
 For $\delta_1$ as chosen above so that \eqref{eq-diskdelta1} holds, set $\e_1' = \e_f(\delta_1)$.

 \begin{lemma}\label{lem-locallifting}
Let $f \colon \fM \to M$ be an $\eF$-map. Let  $x_0 \in \fM$ and suppose that  
$B_{\fM}(x_0,\eU) \subset U_i$ for some   $1 \leq i \leq \nu$. Then for  $z \in B_{\fM}(x_0,\e_1')$ and $y \in D_{\F}(x , \lambda_2) \subset \cP_i(x)$, we have $\ds   f^{-1}(f(y)) \cap \cP_i(z) \ne \emptyset$.
\end{lemma}
 \proof
 
Set $w_0 = f(x_0) \in M$, then  by \eqref{eq-unifcont} we have that   
  $B_{\fM}(x_0,\e_1') \subset f^{-1}(B_{M}(w_0, \delta_1)) \subset U_i$ so that  $f(z) \in B_{M}(w_0, \delta_1)$.
  Moreover,  by Proposition~\ref{prop-coverage}, we have that $B_{M}(w_0, \delta_1) \subset  f(D_{\F}(x_0, \lF/2 ))$.  
  Then choose $x_z \in D_{\F}(x_0, \lF/2 )$ with $f(x_z) = f(z)$.
     
Let  $\gamma_y \colon [0,1] \to \cP_i(x_0)$   be the geodesic  path in  $\cP_i(x_0)$ with  $\gamma_y(0) = x_z$ and $\gamma_y(1) = y$. 

  Let $0 \leq s_* \leq 1$ be the largest value such that 
 $$ f^{-1}(f((\gamma_y(s))) \cap \cP_i(z) \ne \emptyset \quad {\rm for ~ all} \quad  0 \leq s \leq s_* \ .$$ 
 We claim that $s_* =1$. 
Suppose that $s_* < 1$, then we show this yields a contradiction.

Set  $x_* =  \gamma_y(s_*)$ and $w_* = f(x_*)$. Then  there exists  $z_* \in  f^{-1}(w_*) \cap \cP_i(z)$ by the definition of $s_*$.
Note that  $x_* \in  B_{\fM}(x_0,\eF) \subset   B_{\fM}(x_0,\eU/4)$ by the choice of $\lambda_2$ and the fact that $D_{\F}(x_0 , \lambda_2)$ is strongly  convex.
By the choice of $f$ and $\eF$ in \eqref{eq-eF},  we have $\dM(x_* , z_*) < \eF < \eU/4$.
Thus  $\dM(x_0 , z_*) < \eU/2$ and hence $B_{\fM}(z_*, \eF) \subset B_{\fM}(x_0 ,\eU) \subset U_i$.
It then follows from the choice of $\lambda_2$ and the above observations that  
 \begin{equation}\label{eq-disk0}
D_{\F}(z_* , \lambda_2) \subset  B_{\fM}(z_*,\eF)  \subset U_i \ .
\end{equation}

The value of $\eF > 0$ is less than or equal to the choice $\e_1/2$ for this constant in the proof of   Proposition~\ref{prop-coverage} in Appendix~\ref{appendix-A}, and thus 
$B_M(w_*,\delta_1)  \subset f(D_{\F}(z_* , \lF/2))$.  
  The assumption that $s_* <1$ implies that for $s_* \leq s < 1$ sufficiently small so that  $f(\gamma_z(s)) \in B_M(w_*,\delta_1)$, we have that 
    $$  \cP_i(z) \cap f^{-1}(f(\gamma_z(s))) = D_{\F}(z_* , \lF/2) \cap f^{-1}(f(\gamma_z(s))) \ne \emptyset$$
  which contradicts the choice of $s_*$. 
   \endproof

  We next extend the conclusion of Lemma~\ref{lem-locallifting} from paths contained in a coordinate chart, to leafwise paths defined by a plaque chain of arbitrary length.

  Let $\cI = (i_0, i_1, \ldots, i_{\alpha})$ be an admissible chain   with associated holonomy map $h_{\cI} \in \cGF^*$ and $w_0 \in \fD(h_{\cI})$. 
  As in    Section~\ref{subsec-admissible}, we associate to the pair $(\cI, w_0)$ the plaque chain 
$\cP_{\cI}(w_0) = \{\cP_{i_0}(w_0), \cP_{i_1}(w_1), \ldots , \cP_{i_{\alpha}}(w_{\alpha}) \}$ given in \eqref {eq-plaquechain}.

Next   introduce a plaque chain $\cJ = (j_0, j_1, \ldots, j_{\beta})$ which is a refinement of $\cI$ at $w_0$ and is chosen with respect to the leaf distance constant  $\lambda_2 > 0$ which was defined so that the inclusions in \eqref{eq-diameter} hold. 
By   Proposition~\ref{prop-copc}, its associated holonomy map $h_{\cJ}$ at $w_0$ agrees with the holonomy map   $h_{\cI}$ at $w_0$ on their common domains.

Let   $\gamma \colon [0,\alpha] \to L_{x_0} \subset \fM$ be the leafwise piecewise geodesic  associated to the plaque chain $\cP_{\cI}(x_0)$. That is, $\gamma \colon [0,\alpha] \to L_{x_0}$ is the concatenation of geodesic segments $\{\gamma_{\ell} \mid 0 \leq \ell \leq \alpha - 1\}$ in the plaques of the covering $\cU$, where 
$\gamma_{\ell} \colon [\ell  , \ell +1] \to \cP_{i_{\ell}}(w_{\ell})$ satisfies 
$$\gamma_{\ell}(\ell) = x_{\ell} = \tau_{i_{\ell}}(w_{\ell}) \in \cT_{i_{\ell}} \quad , \quad \gamma_{\ell}(\ell +1) = x_{\ell +1} = \tau_{i_{\ell +1}}(w_{i_{\ell +1}})  \in \cT_{i_{\ell +1}} \ .$$

 Introduce a subdivision of the interval $[0,\alpha]$, given by $0 = s_0 <
s_1 < s_2 < \cdots < s_{\beta} = \alpha$, where there is an increasing
subsequence $\{\ell \mid 0 \leq \ell \leq \alpha = s_\beta \}$. For
notational convenience, set $s_{-1} = s_0 = 0$ and $s_{\beta+1} =
s_{\beta} = \alpha$.
Then set $\xi_{\ell} = \gamma(s_{\ell})$ for $-1 \leq \ell \leq \beta+1$, and we
choose the subdivision so ÊÊÊthat Êfor each $0 \leq \ell \leq \beta$,
\begin{equation}\label{eq-pathdistances}
\dF(\gamma(s), \xi_{\ell}) < \lambda_2 ~{\rm for}~ s_{\ell-1} \leq s \leq s_{\ell+1} \ .
\end{equation}

  For each $0 \leq \ell \leq \beta$, choose an index $1 \leq j_{\ell} \leq \nu$ so that $ B_{\fM}(\xi_{\ell}, \eU) \subset U_{j_{\ell}}$. 
   It then follows by the choice of $\eF$, $\lambda_2$ and \eqref{eq-diameter} that  for each $0 \leq \ell \leq \beta$, we have 
\begin{equation}\label{eq-disk2}
\gamma(s) \in  D_{\F}(\xi_{\ell} , \lambda_2) \subset  B_{\fM}(\xi_{\ell},\eF) \subset B_{\fM}(\xi_{\ell}, \eU/4) \quad {\rm for ~ all} ~ s_{\ell-1} \leq s \leq s_{\ell+1} \ .
\end{equation}
Moreover, $\dF(\xi_{\ell}, \xi_{\ell +1}) < \lambda_2$ implies that  
\begin{equation}\label{eq-palquechain2}
\xi_{\ell+1} \in   D_{\F}(\xi_{\ell} , \lambda_2) \cap D_{\F}(\xi_{\ell +1} , \lambda_2) \subset \cP_{j_{\ell}}(\xi_{\ell}) \cap \cP_{j_{\ell +1}}(\xi_{\ell +1}) \ .
\end{equation}
 
Thus     $\cJ  = (j_0, j_1, \ldots, j_{\beta})$ is an admissible sequence, and  
 $\ds \cP_{\cJ}(w_0) = \{\cP_{j_0}(\xi_0), \cP_{j_1}(\xi_1), \ldots , \cP_{j_{\beta}}(\xi_{\beta}) \}$ defines 
  a holonomy map   $h_{\cJ}$    at $w_0$.
 
 Now let $\eF > 0$ be as above,   $\xi_0 = x_0  \in \cT_{j_0}$ for the plaque chain $\cJ$ as chosen above, and 
 suppose that  $\cP_{j_0}(z_0)  \cap f^{-1}(f(\xi_0)) \ne \emptyset$ for some $z_0 \in \fT_{j_0}$.   
 Then $\dF(\xi_0, \xi_1) < \lambda_2$  by \eqref{eq-pathdistances}, so by \eqref{eq-palquechain2} we have 
$\ds \xi_1 \in  \cP_{j_0}(\xi_0) \cap \cP_{j_1}(\xi_1)$. 
Hence by  Lemma~\ref{lem-locallifting} there exists  $\ds  z_1' \in  f^{-1}(f(\xi_1)) \cap \cP_{j_0}(z_0)$.

Note that      $\dM(\xi_1, z_1') \leq \eF \leq \eU/4$, so $z_1' \in B_{\fM}(\xi_1,\eU/4) \subset U_{j_1}$.
  Thus, there exists   $z_1 \in  \cT_{j_1}$   such that $z_1' \in \cP_{j_1}(z_1)$ and hence $\cP_{j_0}(z_0) \cap    \cP_{j_1}(z_1) \ne \emptyset$.

We now repeat the application of  Lemma~\ref{lem-locallifting} to the new basepoint $\xi_1$, and then continue recursively to obtain a sequence of points 
$\ds \{ z_{\ell} \in   \cT_{j_{\ell}}  \mid 0 \leq \ell \leq \beta \}$ such that for $0 < \ell \leq \beta$ we have:
\begin{itemize}
 \item $\cP_{j_{\ell-1}}(z_{\ell-1}) \cap    \cP_{j_{\ell}}(z_{\ell}) \ne \emptyset$ , 
\item $\ds  z_{\ell}' \in  f^{-1}(f(\xi_{\ell})) \cap \cP_{j_{\ell}}(z_{\ell})$.
\end{itemize}

Recall that $\pi_{i} \colon \oU_{i} \to \fT_{i}$ for $1 \leq i \leq \nu$ is the transverse projection to the model space $\fT_i$.
Then the above shows that for $w_0 = \pi_{j_0}(z_0)$ and $w_{\beta} = \pi_{j_{\beta}}(z_{\beta})$ we have $w_0 \in   \fD(h_{\cJ})$ and $h_{\cJ}(w_0) = w_{\beta}$. 

  \subsection{Proof of Theorem~\ref{thm-main1}}\label{subsec-proof}
We can now complete the proof of Theorem~\ref{thm-main1}.  We have assumed that $\e > 0$ is given, and $\eF > 0$ is defined as in \eqref{eq-eF}. 
Then choose an $\eF$-map $f$ as in Section~\ref{subsec-contest}.
Let $h_{\cI} \in \cGF^*$ be as in Section~\ref{subsec-extensions}, and $\cP_{\cJ}$ the path chain constructed above from $\cI$.

For each $1 \leq i \leq \nu$ the transversal map $\tau_i \colon \fT_i \to \cT_i$ is a homeomorphism of compact spaces, 
and the metric $\dTi$ on the subspace $\cT_i \subset \oU_i \subset \fM$ was defined in Section~\ref{subsec-fs}  as the restriction of   $\dM$.

Recall that $\e_1' = \e_f(\delta_1)$ was defined in Section~\ref{subsec-extensions} and used in the hypothesis of Lemma~\ref{lem-locallifting}. 
By the uniform continuity of the maps $\tau_i$, there exists  $\delta > 0$   such that for all $1 \leq i \leq \nu$ and $w \in \fT_i$,
\begin{equation}
B_{\fT_i}(w, \delta) \subset \tau_i^{-1}(B_{\cT_i}(\tau_i(w), \e_1' ))  \ .
\end{equation}
It thus follows from the above results that $h_{\cJ}(B_{\fT_{j_0}}(w, \delta)) \subset B_{\fT_{j_{\beta}}}(h_{\cJ}(w), \e)$, as was to be shown.

\appendix

\section{Local surjectivity for $\e$-maps}\label{appendix-A}

In this appendix, we give a technical result  concerning $\e$-maps.

\begin{prop}\label{Aprop-coverage}
Let $\fM$ be a matchbox manifold with leafwise Riemannian metric on $\F$. Then 
there exists $\eF > 0$  such that,    if $f \colon \fM \to M$ is an $\eF$-map to a compact manifold $M$, then there exists $\delta_1 > 0$ such that for      $x_0 \in \fM$ with $w_0 = f(x_0)$,  we have $D_M(f(x_0),\delta_1) \subset f(D_{\F}(x_0, \lF/2 ))$.
\end{prop}
\proof
   We use the notations of Section~\ref{subsec-metric} above.
   
Choose $\e_0 > 0 $ so that Ê$\rho(\dF,\dM, \e_0) \leq \lF/2$, as
defined by \eqref{eq-rhoFM} to be the maximal radius of a leafwise disk about $x$ 
contained in a closed disk of radius $\e_0$ in $\fM$, for \emph{all} $x \in \fM$. 

Set $\lambda_0 = \lambda_{\F}(\e_0)$, which is
defined in Section \ref{subsec-moremetric} to be the largest
radius $\lambda \leq \lF$ such that the ÊÊdisk Ê$D_{\F}(x, \lambda)$ is
contained in $B_{\fM}(x,\e_0)$, for \emph{all} $x \in \fM$.
Then by the definition of $\rho(\dF,\dM, \e_0)$ and the inclusions         \eqref{eq-metriccont},   for all $x \in \fM$ we have the inclusions
\begin{equation}
D_{\F}(x, \lambda_0 )   \subset     D_{\F}(\dM, x, \e_0) \subset  D_{\F}(x, \lF/2)  \cap    B_{\fM}(x, \e_0) \ .
\end{equation}

Next, choose $\e_1 > 0$ so that  $\rho( \dF, \dM, \e_1) \leq \lambda_0/10$, and let $\lambda_1 = \lF(\e_1)$.
Then  for all $x \in \fM$,  
\begin{equation}\label{Aeq-disklambda0}
D_{\F}(x, \lambda_1)   \subset  D_{\F}(\dM, x, \e_1)  \subset D_{\F}(x, \lambda_0/10)  \cap    B_{\fM}(x, \e_1) \ .
\end{equation}

 Set $ \eF = \e_1/2$. This constant is chosen so that  the result \cite[Section~1, Th\'eor\`eme]{Eilenberg1938} by Eilenberg holds uniformly for strongly convex   compact subsets 
of  $D_{\F}(x, \lF) \subset L_x$, as will be shown below.

 Let $f \colon \fM \to M$ be  an $\eF$-map, which is   onto the compact manifold $M$.  
Let $d_M$ be a metric on $M$. 
For   $w \in M$ and $\delta > 0$, let   $B_{M}(w, \delta) =\{w' \in M \mid d_M(w',w) < \delta\}$  denote the open disk in $M$ of radius $\delta$, and     $D_{M}(w, \delta) =\{w' \in M \mid d_M(w',w) \leq \delta\}$   denote the closed disk in $M$ of radius $\delta$. 

Assume that $d_M$ is chosen so that there exists 
a constant    $\delta_M >0$  such that for all $w \in M$ and $0 < \delta \leq \delta_M$, the   disk $D_{M}(w, \delta)$ is homeomorphic to a disk in $\mR^n$. For example, if $M$ is a Riemannian manifold, then   let $\delta_M > 0$ be   such that each   disk  $D_{M}(w, \delta_M)$ is strongly convex. 

Since $\fM$ is compact, there exists a   uniform modulus of continuity function   $\e_f(\delta) > 0$ for $f$:  for $\delta > 0$, the constant $\e_f(\delta)$ is the largest value  such that 
\begin{equation}\label{Aeq-unifcont}
B_{\fM}(x,\e_f(\delta)) \subset f^{-1}(B_{M}(f(x), \delta)) ~ {\rm for ~ all}~  x \in \fM\ . 
\end{equation}
Let     $\delta_f^* > 0$ be the largest radius such that $\e_f(\delta) \leq \e_{\F}$ for all $0 < \delta \leq \delta_f^*$.

Set $\lambda_{\F,f}(\delta) = \lambda_{\F}(\e_f(\delta))$, which is well-defined for $0 < \delta \leq \delta_f^*$. 
 Combining \eqref{eq-metriccont} and \eqref{Aeq-unifcont} we obtain a leafwise modulus of continuity for $f$:  
 \begin{equation}\label{Aeq-unifleafcont}
  D_{\F}(x, \lambda_{\F,f}(\delta) )   \subset  B_{\fM}(x,\e_f(\delta)) \subset f^{-1}(B_{M}(f(x), \delta)) \quad {\rm for ~ all}~ 0 < \delta \leq \delta_f^* ~, ~ x \in \fM\ .   
\end{equation}

 As $f$ is an $\eF$-map, we have that  $f^{-1}(f(x)) \subset B_{\fM}(x,\eF)$ for all $x \in \fM$. 
As $M$ is compact and $f^{-1}(f(x))$ is a compact set with diameter at most $\eF$ for all $x \in \fM$,  there exists  $0 < \delta_1   \leq \delta_M/10$ so  that  for all $x \in \fM$, we have 
 \begin{equation}\label{Aeq-diskdelta1}
 f^{-1}(D_M(f(x),\delta_1)) \subset B_{\fM}(x,2\eF) = B_{\fM}(x,\e_1) \ .
 \end{equation}

Let $\lambda_2 =  \lambda_{\F,f}(\delta_1)$ so that for all $x \in \fM$,  by \eqref{Aeq-unifleafcont},  \eqref{Aeq-disklambda0},  and the choice of $\delta_1$ we have
  \begin{equation}\label{Aeq-diameter}
D_{\F}(x, \lambda_2 )   \subset f^{-1}(D_M(f(x),\delta_1)) \cap D_{\F}(x,\lF)   \subset D_{\fM}(x,\e_1) \cap D_{\F}(x,\lF)   \subset D_{\F}(x,\lambda_0/10) \ ,
\end{equation}
  \begin{equation}\label{Aeq-diameter2}
  f(D_{\F}(x, \lambda_2 )) \subset D_M(f(x), \delta_1)  \subset D_M(f(x), \delta_M/10) \ .
  \end{equation}

 Set $\lambda_3 = 5 \lambda_2$ so that by \eqref{Aeq-diameter} we have  $\lambda_3 \leq \lF/2$.

\begin{lemma}\label{lem-inclusions}
 Let $x \in \fM$ and   $w = f(x) \in M$, then  $f(D_{\F}(x ,   \lambda_3)) \subset D_M(w, \delta_M)$.
\end{lemma}
\proof
Let $y \in D_{\F}(x, \lambda_3) \subset D_{\F}(x, \lF/2)$ and let $\sigma \colon [0,1] \to D_{\F}(x,\lambda_3)$ be the unique geodesic segment with $\sigma(0) = x$ and $\sigma(1) = y$. 
Set $x_i = \sigma(i/5)$ for $0 \leq i \leq 5$, then $x_{5} = y$. Note that as $\sigma$ is a geodesic, we have 
$$\dF(x_i, x_{i+1}) = \dF(x ,y)/5 \leq \lambda_3/5 = \lambda_2 ~ , ~ 0 \leq i < 5 \ .$$
Set $w_i = f(x_i)$ for $0 \leq i \leq 5$. 

Then for each $0 \leq i \leq 5$, by \eqref{Aeq-diameter2} we have 
$\ds D_{\F}(x_i, \lambda_2) \subset f^{-1}(D_M(w_i,\delta_1))$, 
so that the collection $\{D_M(w_i,\delta_1) \mid 0 \leq i \leq 5\}$ is a covering of the image of $\sigma$. 
Thus, $d_M(w_0 , w_{5}) \leq 10 \delta_1 \leq \delta_M$, hence $d_M(w, f(y)) \leq \delta_M$,  as was to be shown.
\endproof

 For    $x \in \fM$ and $0 < \lambda \leq \lF$,  introduce the following leafwise sets: 
 \begin{eqnarray*}
S_{\F}(x, \lambda) & = &  \{y \in D_{\F}(x, \lF) \mid   \dF(y,x) = \lambda\}   \\
D_{\F}^*(x, \lambda) & = &  \{y \in   D_{\F}(x, \lF) \mid 0 < \dF(y,x) \leq \lambda\} \ .
\end{eqnarray*}
Then we have $S_{\F}(x, \lambda) \subset D_{\F}^*(x, \lambda) \subset D_{\F}(x, \lambda)\subset B_{\F}(x, \lF)$.

Now let  $x_0 \in \fM$ with $w_0 = f(x_0)$. 
 We claim that $B_M(w_0,\delta_1) \subset f(D_{\F}(x_0, \lF/2))$. 
Suppose not, then we show this yields a contraction.

 Let $w_2 \in B_M(w_0,\delta_1) $ but $w_2 \not\in f(D_{\F}(x_0, \lF/2))$. Then there exists $0 < \delta_2 <  \delta_1$ such that 
 \begin{equation}\label{Aeq-delta2}
B_M(w_2,\delta_2)  \subset B_M(w_0,\delta_1) \quad , \quad B_M(w_2,\delta_2) \cap f(D_{\F}(x_0, \lF/2)) = \emptyset \ .
\end{equation}
We consider the maps on \v{C}ech cohomology induced by $f$ to obtain the contradiction.

For $0 < \lambda <    \lF/2$, introduce the collections of open sets in $L_{x_0}$:
\begin{eqnarray}
\cL(\lambda) & = &   \{B_{\F}(x, \lambda)   \mid x \in  D_{\F}(x_0, \lF/2) \}  \\
\cL^*(\lambda) & = &  \{ B_{\F}(x, \lambda) \mid x \in  D_{\F}^*(x_0, \lF/2) ~ {\rm such ~ that} ~ B_{\F}(x, \lambda) \subset D_{\F}^*(x_0, \lF) \}    \ .
\end{eqnarray}

For each $x \in \fM$, the disk $B_{\F}(x, \lambda) \subset L_x$  is strongly convex,   thus $\cL(\lambda)$ is a good covering of $D_{\F}(x_0,\lF/2)$  in the sense of \v{C}ech theory.
Let $\|\cL(\lambda)\|$ denote the simplicial space which is the geometric realization of the collection $\cL(\lambda)$.

For $0 < \lambda < \lambda' < \lF/2$, each open disk $B_{\F}(x, \lambda) \in  \cL(\lambda)$ is contained in the disk $B_{\F}(x, \lambda') \in  \cL(\lambda')$ which induces a   map between  the realizations of their nerve complexes, $\iota \colon  \|\cL(\lambda)\| \to  \|\cL(\lambda')\|$, which is a homotopy equivalence as all the sets in the cover are strongly convex.
Similarly, the induced map on the nerve complex induces  a homotopy equivalence $\iota \colon  \|\cL^*(\lambda)\| \to  \|\cL^*(\lambda')\|$.

 For $0 < \delta < \delta_2$, where $\delta_2$ was chosen so that \eqref{Aeq-delta2} holds, introduce the collections of open sets in $M$:  
  \begin{eqnarray}
\cM_{\delta_2}(\delta) & = &    \{B_M(w, \delta)  \mid w \in   D_M(w_0, \delta_M) - B_M(w_2,\delta_2) \} \\
\cM_f(\delta) & = &    \{B_M(w, \delta)    \mid w \in f(D_{\F}(x_0, \lambda_3)) \} \\
\cM_{S}(\delta) & = &  \{ B_M(w, \delta)    \mid w \in  f(S_{\F}(x_0, \lambda_3) )   \}    \ .
\end{eqnarray}

Recall that by Lemma~\ref{lem-inclusions}, $f(S_{\F}(x, \lambda_3)) \subset D_M(w, \delta_M)$, so   by the choice of $\delta_2$ we have inclusions 
$\cM_{S}(\delta) \subset \cM_f(\delta) \subset \cM^*(\delta)$, and so obtain maps of their simplicial realizations
\begin{equation}\label{Aeq-compositions0}
 \|\cM_{S}(\delta)\| \longrightarrow \|\cM_f(\delta) \| \longrightarrow \|\cM_{\delta_2}(\delta)\|  \ .
\end{equation}

As the disk $D_{\F}(x_0, \lF)$ is strongly convex, there is a natural map $\cR_{\lambda} \colon \|\cL(\lambda)\| \to D_{\F}(x_0, \lF/2)$ which maps a simplex in   the realization $\|\cL(\lambda)\|$    to the geodesic simplex in $D_{\F}(x_0, \lF/2)$ spanned by its vertices. Then $\cR_{\lambda}$ induces   isomorphisms
\begin{eqnarray}
\cR_{\lambda}^*   \colon    \{0\} & \cong &  H^n(D_{\F}(x_0,\lF/2) ; \mZ)   \to   H^n(\|\cL(\lambda)\|; \mZ)  \label{Aeq-compositions1} \\
\cR_{\lambda}^*   \colon    \mZ & \cong &  H^{n-1}(D_{\F}^*(x_0, \lF/2) ; \mZ)   \to   H^{n-1}(\|\cL^*(\lambda)\|; \mZ) \ . \label{Aeq-compositions2}
\end{eqnarray}
   Similarly, for $0 < \delta < \delta_2$,  there is a continuous map $\cS_{\delta_2} \colon \|\cM_{\delta_2}(\delta)\| \to D_M(w, \delta_M) - B_M(w_2,\delta_2)$.
   
  Let  $0 < \delta_3 \leq \delta_2$ be sufficiently small so that we have  an inclusion map 
\begin{equation}\label{Aeq-compositions3}
\cS_{\delta_2}^* \colon \mZ \cong H^{n-1}(\{ D_M(w_0, \delta_M) - B_M(w_2,\delta_2) \} ; \mZ) \to H^{n-1}(\|\cM_{\delta_2}(\delta_3) \| ;\mZ) \ .
\end{equation}

We next consider the maps induced by $f$ on  the cohomology groups in   \eqref{Aeq-compositions1}, \eqref{Aeq-compositions2} and \eqref{Aeq-compositions3}.

Let   $\lambda_4 =  \lambda_{\F,f}(\delta_3)$ so that  for all $x \in \fM$, we have  the inclusion  $ f(D_{\F}(x, \lambda_4)) \subset D_M(f(x), \delta_3)$.
 Then $f$ induces an inclusion map
 $\cU_f \colon \cL(\lambda_4) \to \cM_{\delta_2}(\delta_3)$, which induces a map of their realizations
 \begin{equation}
\| \cU_f\| \colon \| \cL(\lambda_4)\|  \to \| \cM_{\delta_2}(\delta_3) \| \ .
\end{equation}

We have that $0 < \delta_3 \leq \delta_2 <  \delta_1$ so that by  \eqref{Aeq-diskdelta1} and \eqref{Aeq-diameter}, for all $x \in \fM$  we have  
 \begin{equation}\label{Aeq-lifting}
  f^{-1}(D_M(f(x),\delta_3)) \cap D_{\F}(x,\lF)     \subset D_{\F}(x,\lambda_0/10) \ .
\end{equation}
 For each $B_M(w, \delta_3) \in \cM_f(\delta_3)$ choose $x \in f^{-1}(w) \cap D_{\F}(x_0,\lF)$, then the inclusion \eqref{Aeq-lifting} holds. Thus, 
  $f^{-1}$ induces an inclusion map
 $\cV_f \colon  \cM_{\delta_2}(\delta_3) \to \cL(\lambda_0/10)$, which in turn induces a map between their realizations
 \begin{equation}
\| \cV_f\| \colon  \| \cM_{\delta_2}(\delta_3) \| \to \| \cL(\lambda_0/10)\|  \ .
\end{equation}
  Note that the composition $\| \cV_f\| \circ \| \cU_f\| \colon \| \cL(\lambda_4)\| \to \| \cL(\lambda_0/10)\|$ is a homotopy equivalence, as   the sets in $\cL(\lambda_0/10)$ are strongly convex.

Next, for $0 < \lambda \leq \lambda_0/10$,  introduce the   collection  of open balls centered at points of $S_{\F}(x_0, \lambda_3)$, 
\begin{equation}\label{Aeq-sphere}
\cL_{S}(\lambda)   =    \{ B_{\F}(x, \lambda) \mid x \in  S_{\F}(x_0, \lambda_3) ~ {\rm such ~ that} ~ B_{\F}(x, \lambda) \subset D_{\F}^*(x_0, \lF) \} \ .
\end{equation}
As above, the map $f$ induces  maps of geometric realizations 
 \begin{eqnarray}
\| \cU_f\| & \colon & \| \cL_{S}(\lambda_4)\|  \to \| \cM_{S}(\delta_3) \|   \\
\| \cV_f\| & \colon &  \| \cM_{S}(\delta_3) \| \to \| \cL_{S}(\lambda_0/10)\|  \ ,
\end{eqnarray}
and the composition $\| \cV_f\| \circ \| \cU_f\| \colon \| \cL_{S}(\lambda_4)\| \to \| \cL_{S}(\lambda_0/10)\|$ is a homotopy equivalence.

The space $\| \cL_{S}(\lambda_0/10)\|$ is the simplicial realization of a good covering of $S_{\F}(x_0, \lambda_3)$, so we have 
$\ds H^{n-1}( \| \cL_{S}(\lambda_0/10)\|, \mZ) \cong \mZ$.  Let $\omega \in H^{n-1}( \| \cL_{S}(\lambda_0/10)\|, \mZ)$ be a choice of a generator.

Note that the compact sets $f(S_{\F}(x, \lambda))$ in  $M$   limit to the point $w = f(x)$ as $\lambda \to 0$. It follows that there is a non-trivial class 
 $\| \cV_f\|^*(\omega) \in Image \left\{  \cS_{\delta_2}^*  \right\}$ where $\cS_{\delta_2}^*$ is given in \eqref{Aeq-compositions3}. 

On the other hand, $\| \cV_f\|^*$ factors through the map 
$$ \| \cV_f\|^* \colon   H^{n-1}(\| \cL(\lambda_0/10)\| ; \mZ) \to H^{n-1}(\| \cM_{\delta_2}(\delta_3) \| , \mZ) \ .$$
The group $H^{n-1}(\| \cL(\lambda_0/10)\| ; \mZ) \cong \{0\}$, so that $\| \cV_f\|^*(\omega) = 0$, which is a contradiction. 
 \endproof


\end{document}